\documentclass[11pt]{article}

\usepackage{amsfonts}
\usepackage{amsmath}
\usepackage{amssymb}
\usepackage{calc}
\usepackage{color}
\usepackage{setspace}

\newtheorem{theorem}{Theorem}[section]
\newtheorem{lemma}[theorem]{Lemma}

\newenvironment{changemargin}[2]{%
  \begin{list}{}{%
    \setlength{\topsep}{0pt}%
    \setlength{\leftmargin}{#1}%
    \setlength{\rightmargin}{#2}%
    \setlength{\listparindent}{\parindent}%
    \setlength{\itemindent}{\parindent}%
    \setlength{\parsep}{\parskip}%
  }%
  \item[]}{\end{list}}

\newcommand{\qed}{\hspace{1cm} \ensuremath{\Box}}

\makeatletter
\renewcommand{\@cite}[1]{#1}
\makeatother

\title{Manin's conjecture on a nonsingular \\ quartic del Pezzo surface}
\author{Fok-Shuen Leung}
\date{}

\numberwithin{equation}{section}

\begin{document}

\maketitle

\begin{changemargin}{0.5in}{0.5in}

\noindent \small{Given a nonsingular quartic del Pezzo surface,
Manin's conjecture predicts the density of rational points on the
open subset of the surface formed by deleting the lines. We prove
that this prediction is of the correct order of magnitude for a
particular surface.}

\end{changemargin}

\medskip\medskip

\section{Introduction}

\medskip\medskip

Let $V$ be the nonsingular del Pezzo surface of degree four defined
by the zero locus of the equations
\[
           \begin{array}{l}
              0 = x_1x_2- x_3x_4,\\
              0 = x_1^{2}+x_2^{2}+x_3^{2}-x_4^{2}-2x_5^{2}.
           \end{array}
\]
Let $U\subseteq V$ be formed by deleting the lines from $V$. Given a
rational point
$\textbf{x}=[x_1,\ldots,x_5]\in\mathbb{P}^4(\mathbb{Q})$ with
$x_1,\ldots,x_5\in\mathbb{Z}$ and
$\gcd\left(x_1,\ldots,x_5\right)=1$, we define the height of
$\textbf{x}$ to be
$\|\textbf{x}\|=\max\left(|x_1|,\ldots,|x_5|\right)$. Given $B\ge1$,
the density of rational points on $V$ is specified by the
cardinality
\[
N_U(B)=\#\left\{\textbf{x}\in U\cap\mathbb{P}^4(\mathbb{Q}) :
\|\textbf{x}\|\le B\right\}.
\]
Manin's conjecture, proposed in [\cite{Franke/Manin/Tschinkel 1989}]
for Fano varieties in general, predicts in this case that
\[
N_{U}(B)=c_V B \left(\log B\right)^{\rho-1}\left(1+o(1)\right)
\]
as $B\rightarrow\infty$, where $c_V$ is a positive constant and
$\rho$ is the rank of the Picard group of $V$. Our principal result
is the following:

\begin{theorem}\label{1 thm main}
$B\left(\log B\right)^{\rho-1}\ll N_{U}(B)\ll B\left(\log
B\right)^{\rho-1}$.
\end{theorem}

\medskip\medskip

An overview of progress in proving Manin's conjecture for del Pezzo
surfaces can be found in [\cite{Browning 2006}]. In general,
singular del Pezzo surface of low degree have proven more tractable
than their nonsingular counterparts. For nonsingular quartic
surfaces, the best result until now is due to Salberger, who proved
$N_{U}(B)\ll B\left(\log B\right)^{1+\varepsilon}$ for any
$\varepsilon>0$, provided $V$ contains a rational conic; this work
was presented at the 2001 Budapest conference \textit{Higher
dimensional varieties and rational points}. Our result refines
Salberger's.

\medskip\medskip

Both bounds comprising Theorem \ref{1 thm main} involve fibering $V$
into a family of conics; this allows us to reduce the problem of
estimating $N_{U}(B)$ to the problem of estimating the density of
certain rational points on these conics. The same idea is central to
Salberger's result; our improved bound stems from tighter control on
the uniformity of bounds for rational points on the conics. The
method appears to be applicable in a far more general setting, and
we intend to explore this in a future paper.

\medskip\medskip

\section{The constant $\rho$}

\medskip\medskip

We begin by recounting some geometry of quartic del Pezzo surfaces.
We refer the reader to [\cite{Manin 1974}] for a comprehensive
exposition.

\medskip\medskip

In general, a nonsingular quartic del Pezzo surface $X$ contains 16
lines, each of which intersects exactly five others. Given any
subset of five pairwise skew lines $L_1,\ldots,L_5$, $X$ is
isomorphic to $\mathbb{P}^2$ blown up along five points
$P_1,\ldots,P_5$ in general position such that $L_1,\ldots,L_5$ are
the preimages of those points under the blowup. Moreover, there
exists a unique line $L_0$ intersecting $L_1,\ldots,L_5$; $L_0$ is
the preimage of the unique conic on $\mathbb{P}^2$ through
$P_1,\ldots,P_5$.

\medskip\medskip

Let $K_0,\ldots,K_6$ denote the linear equivalence classes of
$L_0,\ldots,L_6$, respectively, and $K$ denote the class of the
preimage of a line on $\mathbb{P}^2$. Then
\begin{equation}\label{2 eq classes}
K_0\sim 2K-\left(K_1+\cdots+K_5\right).
\end{equation}
The geometric Picard group of $X$
--- that is, the Picard group of $X$ defined over an extension
$E$ of minimal degree over $\mathbb{Q}$ such that all the lines on
$X$ are defined over $E$ --- has a basis $\{K,K_1,\ldots,K_5\}$. The
Picard group of $X$ is that part of the geometric Picard group
invariant under the action of $\mbox{Gal}\left(E/\mathbb{Q}\right)$.

\medskip\medskip

The 16 lines on $V$ have the following parametrizations:

\[
\begin{array}{l l l l}
L_{1}: & \left[a,b,a,b,a\right], & L_{2}: & \left[a,b,a,b,-a\right],
\smallskip
\\
L_{3}: & \left[a,b,-a,-b,a\right], & L_{4}: &
\left[a,b,-a,-b,-a\right],
\smallskip
\\
L_{5}: & \left[a,b,b,a,b\right], & L_{6}: & \left[a,b,b,a,-b\right],
\smallskip
\\
L_{7}: & \left[a,b,-b,-a,b\right], & L_{8}: &
\left[a,b,-b,-a,-b\right],
\smallskip
\\
L_{9}: & \left[a,b,ia,-ib,b\right], & L_{10}: &
\left[a,b,ia,-ib,-b\right],
\smallskip
\\
L_{11}: & \left[a,b,-ia,ib,b\right], & L_{12}: &
\left[a,b,-ia,ib,-b\right],
\smallskip
\\
L_{13}: & \left[a,b,-ib,ia,a\right], & L_{14}: &
\left[a,b,-ib,ia,-a\right],
\smallskip
\\
L_{15}: & \left[a,b,ib,-ia,a\right], & L_{16}: &
\left[a,b,ib,-ia,-a\right].
\smallskip
\end{array}
\]

\noindent Note that all the lines are defined over $\mathbb{Q}(i)$.
Let $K_0,\ldots,K_5$ denote the classes of $L_5$, $L_1$, $L_4$,
$L_6$, $L_9$ and $L_{11}$, respectively. Note that the latter five
lines are pairwise skew, and that they are intersected by $L_5$. Let
$K$ denote the class of the preimage of a line on $\mathbb{P}^2$. In
view of (\ref{2 eq classes}), since $K_0$, $K_1$, $K_2$ $K_3$ and
$K_4+K_5$ are invariant under the action of
$\mbox{Gal}\left(\mathbb{Q}(i)/\mathbb{Q}\right)$, so too is $K$;
and since $\{K, K_1,\ldots,K_5\}$, being a basis, is a linearly
independent set, the set $\{K, K_1, K_2, K_3, K_4+K_5\}$ is also
linearly independent. Therefore the Picard group of $V$ has rank at
least 5. Since not all the lines on $V$ are invariant under the
action of $\mbox{Gal}\left(\mathbb{Q}(i)/\mathbb{Q}\right)$, we
conclude that the Picard group of $V$ has rank exactly 5.

\medskip\medskip

\section{The lower bound}

\medskip

\subsection{Preliminaries}

\medskip\medskip

Let $B>0$ be given and
\[
P=\left\{ (r,s): \mbox{$s$ is even, $\gcd(r,s)=1$ and $1\le r,s\le
B^{1/100}$}\right\}.
\]
Given $(r,s)\in P$, the first quadric of $V$ is satisfied by taking
$x_{1}=rX_{1}$, $x_{2}=sX_{2}$, $x_{3}=sX_{1}$ and $x_{4}=rX_{2}$;
and setting $x_{5}=X_{3}$, the second quadric of $V$ is a ternary
quadric $0=Q_{r,s}(\textbf{X})$, where
\[
Q_{r,s}(\textbf{X})=(r^{2}+s^{2})X_{1}^{2}-(r^{2}-s^{2})X_{2}^{2}-2X_{3}^{2}.
\]
If $\gcd\left(X_1,X_2,X_3\right)=1$, then $\|\textbf{x}\|\le B$ is
implied by the bounds
\begin{equation}\label{3LB/eq X height}
|X_{1}|,|X_{2}|\leq \frac{B}{\mbox{max}(r,s)}, |X_{3}|\leq B.
\end{equation}
Let
\[
N_{r,s}=\# \left\{ \mbox{$\textbf{X}$ : $0=Q_{r,s}(\textbf{X})$,
$\gcd\left(X_1,X_2,X_3\right)=1$ and (\ref{3LB/eq X height}) holds}
\right\},
\]
and let $P_{i}$ denote the set of pairs $(r,s)\in P$ in the dyadic
ranges
\[
2^{i-1}=R_{i}<r\le 2R_{i} = 2^{i}, 2^{i}=S_{i}<s\le 2S_{i} =
2^{i+1}.
\]
(Note that, given $(r,s)\in P_i$ for any $i$, we have $r<s$.) Then
\begin{equation}\label{3LB/eq comp N_UB, N_rs}
N_{U}(B)\gg \sum_{i}\sum_{(r,s)\in P_i}N_{r,s},
\end{equation}
where the $i$ are summed over those values such that the sets
$P_{i}$ are nonempty.

\medskip

\subsection{The cardinality $N_{r,s}$}

\medskip\medskip

Let $(r,s)\in P_i$ be given. We estimate $N_{r,s}$ by parametrizing
a subset of rational points on the quadric $0=Q_{r,s}(\textbf{X})$.

\medskip\medskip

We begin by observing that $[1,1,s]$ is a point on
$0=Q_{r,s}(\textbf{X})$. We fix a nonzero integer constant $c$, and
consider all points on the quadric of the form
\[
\textbf{X}=[c+x+2sy,c+x,cs],
\]
where $(x,y)$ is an integer pair satisfying the coprimality
condition
\begin{equation}\label{3LB/eq (x,y) coprimality}
\gcd(x,2sy)=1.
\end{equation}
Note that distinct pairs $(x,y)$ parametrize distinct points
$\textbf{X}$. We proceed to eliminate the constant $c$. Substituting
$\textbf{X}$ back into $0=Q_{r,s}(\textbf{X})$, we get
\[
0=\left(r^2+s^2\right)\left((x+2sy)^2+2c(x+2sy)\right)-\left(r^2-s^2\right)\left(x^2+2cx\right).
\]
We rearrange this to get
\[
cf_{r,s}(x,y)=-\left(r^2+s^2\right)(x+2sy)^2+\left(r^2-s^2\right)x^2,
\]
where
\[
f_{r,s}(x,y)=2\left(r^2+s^2\right)(x+2sy)-2x\left(r^2-s^2\right).
\]
We simplify $\textbf{X}$ by multiplying each of its components by
$f_{r,s}(x,y)$ and then dividing out by $s^2$, getting
$\textbf{X}=[f_{1,r,s}(x,y),f_{2,r,s}(x,y),f_{3,r,s}(x,y)]$, where
\[
\begin{array}{r c l}
f_{1,r,s}(x,y) &=& x^{2}+4sxy+2\left(r^{2}+s^{2}\right)y^{2},
\smallskip
\\
f_{2,r,s}(x,y) &=& x^{2}-2\left(r^{2}+s^{2}\right)y^{2},
\smallskip
\\
f_{3,r,s}(x,y) &=&
-sx^{2}-2\left(r^{2}+s^{2}\right)xy-2s\left(r^{2}+s^{2}\right)y^{2}.
\end{array}
\]

\medskip\medskip

Now given an integer pair $(x,y)$ satisfying (\ref{3LB/eq (x,y)
coprimality}), the forms $f_{1,r,s}(x,y)$, $f_{2,r,s}(x,y)$ and
$f_{3,r,s}(x,y)$ may have a nontrivial common divisor:

\begin{lemma}\label{3LB/lem gcd(f_1,f_2,f_3) simplification}
Let $(x,y)$ be an integer pair satisfying $(\ref{3LB/eq (x,y)
coprimality})$. Then the greatest common divisor of
$f_{1,r,s}(x,y)$, $f_{2,r,s}(x,y)$ and $f_{3,r,s}(x,y)$ is equal to
\[
\gcd\left(x,r^2+s^2\right)\gcd\left(x+2sy,r^2-s^2\right).
\]
\end{lemma}

\noindent PROOF. Note that
\[
f_{1,r,s}(x,y)+f_{2,r,s}(x,y)=2x(x+2sy).
\]
Now $2$, $x$ and $x+2sy$ are pairwise coprime; hence the greatest
common divisor of $f_{1,r,s}(x,y)$, $f_{2,r,s}(x,y)$ and
$f_{3,r,s}(x,y)$ is equal to the product of the factors
$\gcd\left(2,f_{2,r,s}(x,y),f_{3,r,s}(x,y)\right)$,
$\gcd\left(x,f_{2,r,s}(x,y),f_{3,r,s}(x,y)\right)$ and
$\gcd\left(x+2sy,f_{2,r,s}(x,y),f_{3,r,s}(x,y)\right)$. We denote
these factors $F_1$, $F_2$ and $F_3$, respectively, and simplify
each in turn. For the first, (\ref{3LB/eq (x,y) coprimality})
implies that $x$, hence $f_{2,r,s}(x,y)$, is odd; thus $F_1=1$. For
the second, we again apply (\ref{3LB/eq (x,y) coprimality}), getting
\[
F_2 =
\gcd\left(x,2\left(r^{2}+s^{2}\right)y^{2},2s\left(r^{2}+s^{2}\right)y^{2}\right)
=\gcd\left(x,r^{2}+s^{2}\right).
\]
For the third, note that $f_{2,r,s}(x,y) =
(x+2sy)(x-2sy)-2\left(r^{2}-s^{2}\right)y^2$ and $f_{3,r,s}(x,y) =
-(x+2sy)\left(sx+2r^2y\right)+2s\left(r^{2}-s^{2}\right)y^2$; hence
\[
F_3=
\gcd\left(x+2sy,2\left(r^{2}-s^{2}\right)y^2,2s\left(r^{2}-s^{2}\right)y^2\right)
=\gcd\left(x+2sy,r^{2}-s^{2}\right),
\]
which completes the proof. \qed

\medskip\medskip

Let $\gcd\left(x,r^2+s^2\right)\gcd\left(x+2sy,r^2-s^2\right)=n$.
Then, given a point
$\textbf{X}=[f_{1,r,s}(x,y),f_{2,r,s}(x,y),f_{3,r,s}(x,y)]$, the
bounds (\ref{3LB/eq X height}) are implied by the bounds
\[
\frac{|f_{1,r,s}(x,y)|}{n},\frac{|f_{2,r,s}(x,y)|}{n}\le\frac{B}{s},\frac{|f_{3,r,s}(x,y)|}{n}\le
B,
\]
which are themselves implied by the bounds
\[
1\le x\le X= \left(\frac{Bn}{4s}\right)^{1/2}, |y|\le Y=
\left(\frac{Bn}{16s^3}\right)^{1/2}.
\]
For convenience, we let $z=x+2sy$, which allows us to replace the
above bounds with $1\le x,z\le X$.

\medskip\medskip

We estimate $N_{r,s}$ by indexing the pairs $(x,y)$ contributing to
$N_{r,s}$ according to the greatest common divisor of the components
of $\textbf{X}$. Let
\[
N_{n,r,s}=\#\left\{
\begin{array}{l}
\gcd(x,sz)=1, 2s|x-z, 1\le x,z\le X,
\smallskip
\\
\gcd\left(x,r^2+s^2\right)\gcd\left(x+2sy,r^2-s^2\right)=n
\end{array}
\right\}.
\]
Then
\[
N_{r,s}\ge\displaystyle\sum_{n\ge1}N_{n,r,s}.
\]

\medskip\medskip

The most cumbersome condition on $N_{n,r,s}$ is the last. In order
to keep track of it, we redefine $N_{n,r,s}$ in terms of positive
integer pairs $(a,b)$, where $\gcd\left(x,r^2+s^2\right)=a$,
$\gcd\left(z,r^2-s^2\right)=b$ and $ab=n$. We write $x=au$,
$r^2+s^2=ac$, $z=bv$ and $r^2-s^2=bd$, where
\begin{equation}\label{3LB/eq N_ab condition 1}
\gcd\left(u,c\right)=1\hspace{0.1in}\mbox{and}\hspace{0.1in}\gcd\left(v,d\right)=1.
\end{equation}
The last condition on $N_{n,r,s}$ is implicit in these definitions.
The coprimality condition $\gcd(x,sz)=1$ is implied by
\begin{equation}\label{3LB/eq N_ab condition 2}
\gcd(a,v)=\gcd(u,b)=1=\gcd(u,v)=1=\gcd(u,s)=1;
\end{equation}
the divisibility condition $2s|x-z$ is simply restated
\begin{equation}\label{3LB/eq N_ab condition 3}
2s|au-bv;
\end{equation}
and the bounds $1\le x,z\le X$ are implied by the bounds
\begin{equation}\label{3LB/eq N_ab condition 4}
1\le u\le U=\left(\frac{Bb}{4as}\right)^{1/2}, 1\le v\le
V=\left(\frac{Ba}{4bs}\right)^{1/2}.
\end{equation}
Thus, defining
\[
N_{a,b,r,s}=\#\left\{ (u,v) : \mbox{(\ref{3LB/eq N_ab condition 1}),
(\ref{3LB/eq N_ab condition 2}), (\ref{3LB/eq N_ab condition 3}) and
(\ref{3LB/eq N_ab condition 4}) hold} \right\},
\]
we have
\begin{equation}\label{3LB/eq comp N_UB, N_abrs}
N_{r,s}\ge\sum_{a|r^2+s^2}\sum_{b|r^2-s^2}N_{a,b,r,s}.
\end{equation}

\medskip

\subsection{The cardinality $N_{a,b,r,s}$}

\medskip\medskip

Let $(r,s)\in P_i$, $a|r^2+s^2$ and $b|r^2-s^2$ be given. We
estimate $N_{a,b,r,s}$ by fixing $u$ and then estimating the number
of $v$ such that $(u,v)$ contributes to $N_{a,b,r,s}$.

\medskip\medskip

Given $u$ such that $\gcd(u,s)=1$, let
\[
N_{u,a,b,r,s}=\#\left\{
\begin{array}{l}
\gcd(v,d)=\gcd(v,a)=\gcd(v,u)=1,
\smallskip
\\
2s|au-bv, 1\le v\le V
\end{array} \right\}.
\]
Then
\[
N_{a,b,r,s}=\sum_{u}N_{u,a,b,r,s},
\]
where the sum is taken over a suitable set of $u$. We shall define
this set below.

\medskip\medskip

We use the M\"{o}bius function to pick out the coprimality
conditions on $N_{u,a,b,r,s}$. Let
\[
N'_{u,a,b,r,s}(n_1,n_2,n_3)=\#\left\{ v:
\mbox{$\mbox{lcm}(n_1,n_2,n_3)|v$, $2s|au-bv$, $1\le v\le V$}
\right\}.
\]
Then
\[
N_{u,a,b,r,s}=\sum_{n_{1}|d}\sum_{n_{2}|a}\sum_{n_3|u}\mu(n_{1})\mu(n_{2})\mu(n_3)N'_{u,a,b,r,s}(n_1,n_2,n_3).
\]
Let $n_{1}$, $n_{2}$ and $n_3$ be in the range of summation above.
Then $\gcd(2s,n_1)=1$, since $n_1|r^2-s^2$ and $\gcd(2s,r^2-s^2)=1$;
$\gcd(2s,n_2)=1$, since $n_2|r^2+s^2$ and $\gcd(2s,r^2+s^2)=1$; and
$\gcd(2s,n_3)=1$, since $\gcd(u,s)=1$ and $s$ is even. Moreover,
$\gcd(2s,b)=1$, since $b|r^2-s^2$. Thus
\[
N'_{u,a,b,r,s}(n_1,n_2,n_3)=\frac{V}{2s\cdot\mbox{lcm}(n_1,n_2,n_3)}+O(1),
\]
and
\[
N_{u,a,b,r,s}=\sum_{n_{1}|d}\sum_{n_{2}|a}\sum_{n_3|u}\mu(n_{1})\mu(n_{2})\mu(n_3)
\left(\frac{V}{2s\cdot\mbox{lcm}(n_1,n_2,n_3)}+O(1)\right).
\]
We estimate $N_{a,b,r,s}$ by summing $N_{u,a,b,r,s}$ over the set
\[
P_{a,b,r,s}(n_3)=\left\{ u : \gcd(u,c)=\gcd(u,b)=\gcd(u,s)=1, n_3|u,
1\le u\le U \right\};
\]
that is, $N_{a,b,r,s}$ is equal to
\[
\displaystyle\sum_{n_{1}|d}\sum_{n_{2}|a}\sum_{n_3\le U}\sum_{u\in
P_{a,b,r,s}(n_3)}\mu(n_{1})\mu(n_{2})\mu(n_3)\left(\frac{V}{2s\cdot\mbox{lcm}(n_1,n_2,n_3)}+O(1)\right).
\]
Since the cardinality of $P_{a,b,r,s}(n_3)$ has an upper bound
$U/n_3$, the contribution to $N_{a,b,r,s}$ of the error term above
is of order at most
\[
U\sum_{n_{1}|d}\sum_{n_{2}|a}\sum_{n_3\le U} \frac{1}{n_3}\le
U\left(R_iS_iU\right)^{\varepsilon}
\]
for any $\varepsilon>0$, provided $i$ and $B$ are sufficiently
large; that is,
\[
N_{a,b,r,s}=\frac{V}{2s}\displaystyle\sum_{n_{1}|d}\sum_{n_{2}|a}\sum_{n_3\le
U}\sum_{u\in
P_{a,b,r,s}(n_3)}\frac{\mu(n_{1})\mu(n_{2})\mu(n_3)}{\mbox{lcm}(n_1,n_2,n_3)}
+O\left(U\left(R_iS_iU\right)^{\varepsilon}\right).
\]

\medskip\medskip

We now estimate the cardinality of $P_{a,b,r,s}(n_3)$ more
precisely. As in the case of $N_{u,a,b,r,s}$, we use the M\"{o}bius
function to pick out coprimality conditions on the set. Let
\[
P'_{a,b,r,s}(n_3,m_1,m_2,m_3)=\left\{ u:
\mbox{lcm}(n_3,m_1,m_2,m_3)|u, 1\le u\le U \right\}.
\]
Then
\[
\#P_{a,b,r,s}=\displaystyle\sum_{m_{1}|c}\sum_{m_{2}|b}\sum_{m_3|s}\mu(m_{1})\mu(m_{2})\mu(m_3)\#P'_{a,b,r,s}(n_3,m_1,m_2,m_3).
\]
Now
\[
\#P'_{a,b,r,s}(n_3,m_1,m_2,m_3)=\frac{U}{\mbox{lcm}(n_3,m_1,m_2,m_3)}+O(1).
\]
The contribution to $N_{a,b,r,s}$ of the error term is of order at
most
\[
\displaystyle\frac{V}{s}\displaystyle\sum_{n_{1}|d}\sum_{n_{2}|a}\sum_{n_3\le
U} \sum_{m_{1}|c}\sum_{m_{2}|b}\sum_{m_3|s}
\frac{1}{\mbox{lcm}(n_1,n_2,n_3)} \le
\displaystyle\frac{V}{s}\left(R_iS_iU\right)^{\varepsilon}
\]
for any $\varepsilon>0$, provided $i$ and $B$ are sufficiently
large; that is, $N_{a,b,r,s}$ is equal to
\[
\begin{array}{l}
\displaystyle\frac{UV}{2s}\displaystyle\sum_{n_{1}|d}\sum_{n_{2}|a}\sum_{n_3\le
U}
\sum_{m_{1}|c}\sum_{m_{2}|b}\sum_{m_3|s}\displaystyle\frac{\mu(n_{1})\mu(n_{2})\mu(n_3)\mu(m_{1})\mu(m_{2})\mu(m_3)}
{\mbox{lcm}(n_1,n_2,n_3)\cdot\mbox{lcm}(n_3,m_1,m_2,m_3)}
\medskip
\\
\hspace{0.5in}+O\left(\left(U+\displaystyle\frac{V}{s}\right)\left(R_iS_iU\right)^{\varepsilon}\right).
\end{array}
\]

\medskip\medskip

Finally we estimate the main term above. Let $T_{a,b,r,s}$ denote
this term, and $T'_{a,b,r,s}$ denote $T_{a,b,r,s}$ but with the
difference that, in $T'_{a,b,r,s}$, $n_3$ is summed over all
positive integers rather than over the range $n_3\le U$. Now
$T'_{a,b,r,s}-T_{a,b,r,s}$ is of order at most
\[
\displaystyle\frac{UV}{s}\displaystyle\sum_{n_{1}|d}\sum_{n_{2}|a}\sum_{n_3>U}
\sum_{m_{1}|c}\sum_{m_{2}|b}\sum_{m_3|s}\displaystyle\frac{1}{n_3^2}
\le
\displaystyle\frac{UV}{s}(R_iS_i)^\varepsilon\sum_{n_3>U}\frac{1}{n_3^2}
\le \displaystyle\frac{V}{s}(R_iS_i)^\varepsilon
\]
for any $\varepsilon>0$, provided $i$ and $B$ are sufficiently
large; that is,
\[
N_{a,b,r,s}=T'_{a,b,r,s}+O\left(\left(U+\displaystyle\frac{V}{s}\right)\left(R_iS_iU\right)^{\varepsilon}\right).
\]
In order to estimate $T'_{a,b,r,s}$, we define the condition
\begin{equation}\label{3LB/eq T' divisibility conditions}
n_{1}|d,\hspace{0.1in}n_{2}|a,\hspace{0.1in}m_{1}|c,\hspace{0.1in}m_{2}|b\hspace{0.1in}\mbox{and}\hspace{0.1in}m_3|s,
\end{equation}
and the function $f_{a,b,r,s}(n_1,n_2,n_3,m_1,m_2,m_3)$ to be equal
to
\[
\left\{
\begin{array}{ll}
\displaystyle\frac{\mu(n_1)\mu(n_2)\mu(n_3)\mu(m_1)\mu(m_2)\mu(m_3)}{\mbox{lcm}(n_1,n_2,n_3)\mbox{lcm}(n_3,m_1,m_2,m_3)}
&\mbox{if (\ref{3LB/eq T' divisibility conditions}) holds}
\smallskip
\\
0 & \mbox{otherwise}
\end{array}
\right..
\]
Then
\[
T'_{a,b,r,s}=\displaystyle\frac{UV}{2s}\sum_{\begin{subarray}{c}n_i,m_i\ge1\\\mbox{\scriptsize{for
$1\le i\le3$}}\end{subarray}}f_{a,b,r,s}(n_1,n_2,n_3,m_1,m_2,m_3).
\]
Because $f_{a,b,r,s}$ is multiplicative and we have
\[
\sum_{\begin{subarray}{c}n_i,m_i\ge1\\\mbox{\scriptsize{for $1\le
i\le3$}}\end{subarray}}|f_{a,b,r,s}(n_1,n_2,n_3,m_1,m_2,m_3)|\le\displaystyle\sum_{n_{1}|d}\sum_{n_{2}|a}
\sum_{m_{1}|c}\sum_{m_{2}|b}\sum_{m_3|2s}\sum_{n_3\ge1}
\frac{1}{n_3^2},
\]
which converges, we may write
\[
T'_{a,b,r,s}=\displaystyle\frac{UV}{2s}\prod_{p}f_{p,a,b,r,s},
\]
where the product is taken over all primes $p$, and the local
factors $f_{p,a,b,r,s}$ are defined
\[
f_{p,a,b,r,s} =
\displaystyle\sum_{\begin{subarray}{c}e_i,e'_i\in\{0,1\}\\\mbox{\scriptsize{for
$1\le i\le3$}}
\end{subarray}}f_{a,b,r,s}\left(p^{e_1},p^{e_2},p^{e_3},p^{e'_1},p^{e'_2},p^{e'_3}\right).
\]
We evaluate $f_{p,a,b,r,s}$ directly, in three cases. If $p$ does
not divide any element in the set $\{a,b,c,d,s\}$, then
$f_{p,a,b,r,s}=1-p^{-2}$; if $p$ divides exactly one element in the
set $\{a,b,c,d,s\}$, then $f_{p,a,b,r,s}=1-p^{-1}$; and if $p$
divides exactly two elements in the set $\{a,b,c,d,s\}$ --- that is,
either $p|a$ and $p|c$, or $p|b$ and $p|d$ --- then
$f_{p,a,b,r,s}=\left(1-p^{-1}\right)^2$. Hence
\[
\begin{array}{r c l}
T'_{a,b,r,s} &\ge& \displaystyle\frac{UV}{2s} \prod_{p\nmid
s\Delta_{r,s}}\left(1-\frac{1}{p^2}\right)
\prod_{p|s}\left(1-\frac{1}{p}\right) \prod_{p|\Delta_{r,s}}
\left(1-\frac{1}{p}\right)^2\smallskip
\\
&\gg& \displaystyle\frac{UV}{s}\prod_{p|s}
\left(1-\frac{1}{p}\right)\prod_{p|\Delta_{r,s}}
\left(1-\frac{1}{p}\right)^2,
\end{array}
\]
where $\Delta_{r,s}$ denotes $|r^4-s^4|$, and the relation $\gg$
does not depend on our choice of $a$, $b$, $r$ or $s$. (For the
remainder of this section we assume that all relations $\gg$ are
thus independent.)

\medskip\medskip

Thus
\[
N_{a,b,r,s}\gg\displaystyle\frac{UV}{s}\prod_{p|s}
\left(1-\frac{1}{p}\right)\prod_{p|\Delta_{r,s}}
\left(1-\frac{1}{p}\right)^2+O\left(\left(U+\displaystyle\frac{V}{s}\right)\left(R_iS_iU\right)^{\varepsilon}\right)
\]
for any $\varepsilon>0$, provided $i$ and $B$ are sufficiently
large. We conclude that
\begin{equation}\label{3LB/eq N_abrs}
N_{a,b,r,s} \gg \displaystyle\frac{B}{s^2}\prod_{p|s}
\left(1-\frac{1}{p}\right)\prod_{p|\Delta_{r,s}}
\left(1-\frac{1}{p}\right)^2.
\end{equation}

\medskip

\subsection{The cardinality $N_{U}(B)$}

\medskip\medskip

For convenience we define the multiplicative function
\[
f(n)=\displaystyle\prod_{p|n}\left(1-\frac{1}{p}\right)
\]
for any $n\in\mathbb{N}$, with $f(1)=1$. With this notation, and in
view of the bounds (\ref{3LB/eq comp N_UB, N_rs}), (\ref{3LB/eq comp
N_UB, N_abrs}) and $(\ref{3LB/eq N_abrs})$, we have
\[
\begin{array}{rcl}
N_{U}(B) &\gg& B\displaystyle\sum_{i}\frac{1}{S_i^2}\sum_{(r,s)\in
P_i}\sum_{\begin{subarray}{r}a|r^2+s^2\\b|r^2-s^2\end{subarray}}
f(s)f(\Delta_{r,s})^2
\\
&\ge& B\displaystyle\sum_{i}\frac{1}{S_i^2}\sum_{(r,s)\in P_i}
d\left(\Delta_{r,s}\right)f(s)f(\Delta_{r,s})^2,
\end{array}
\]
where the $i$ are summed over those values such that the $P_i$ are
nonempty. We may restrict the range of summation on the right-hand
side above without invalidating the bound, and it will be useful to
impose the condition that, for any pair $(r,s)$ in that range of
summation, $s$ is not only even but divisible by $6$; that is,
\begin{equation}\label{3LB/eq N_UB 1}
N_{U}(B) \gg
B\displaystyle\sum_{i}\frac{1}{S_i^2}\displaystyle\sum_{\begin{subarray}{c}S_i<s\le2S_i\\6|s\end{subarray}}f(s)
\sum_{\begin{subarray}{c}R_i<r\le2R_i\\\gcd(r,s)=1\end{subarray}}d\left(\Delta_{r,s}\right)f(\Delta_{r,s})^2.
\end{equation}

\medskip\medskip

We estimate the inner sum on the right-hand side of (\ref{3LB/eq
N_UB 1}). Let $s$ be in the range of summation. By the M\"{o}bius
inversion formula, we have
\[
d\left(n\right)f(n)^2=\sum_{m|n}f'(m)
\]
for any $n\in\mathbb{N}$ if, and only if,
\[
f'(n)=\displaystyle\sum_{m|n}\mu\left(\frac{n}{m}\right)d(m)f(m)^2
\]
for any $n\in\mathbb{N}$. Now $f'$ is multiplicative, and given a
prime power $p^e$ with $e\ge1$, we have
\[
f'\left(p^e\right) = \left\{
\begin{array}{l l}
2\left(1-\displaystyle\frac{1}{p}\right)^2-1 &\mbox{if $e=1$}
\smallskip
\\
\left(1-\displaystyle\frac{1}{p}\right)^2 &\mbox{otherwise}
\end{array}
\right.;
\]
that is, $f'\left(p^e\right)>0$ for any $e\in\mathbb{N}$ provided
$p\ge5$. No primes smaller than $5$ divide $\Delta_{r,s}$, since $2$
and $3$ both divide $s$; hence
\[
d\left(\Delta_{r,s}\right)f(\Delta_{r,s})^2=\sum_{m|\Delta_{r,s}}f'(m)\ge\sum_{\begin{subarray}{c}m|\Delta_{r,s}\\m\le
R_i^{1/2}\end{subarray}}f'(m),
\]
where $f'(m)$ is nonnegative over the range of summation. (It will
shortly become clear why we impose a bound on $m$.) Thus
\[
\sum_{\begin{subarray}{c}R_i<r\le2R_i\\\gcd(r,s)=1\end{subarray}}d\left(\Delta_{r,s}\right)f(\Delta_{r,s})^2
\ge
\sum_{\begin{subarray}{c}R_i<r\le2R_i\\\gcd(r,s)=1\end{subarray}}\sum_{\begin{subarray}{c}m|\Delta_{r,s}\\m\le
R_i^{1/2}\end{subarray}}f'(m).
\]

\medskip\medskip

We use the M\"{o}bius function to pick out the coprimality condition
on the right-hand side. As an intermediate step, we define
\[
N_{m,s}=\#\left\{ r: R_i<r\le2R_i, \gcd(r,s)=1, m|\Delta_{r,s}
\right\}.
\]
Then
\[
\sum_{\begin{subarray}{c}R_i<r\le2R_i\\\gcd(r,s)=1\end{subarray}}d\left(\Delta_{r,s}\right)f(\Delta_{r,s})^2
\ge\sum_{\begin{subarray}{c}m\le
R_i^{1/2}\\\gcd(m,s)=1\end{subarray}}f'(m)N_{m,s}.
\]
We impose the condition that $\gcd(m,s)=1$ on the range of summation
on the right-hand side to ensure that the $N_{m,s}$ we sum are
nonzero.

\medskip\medskip

Let $N_{m,s}$ be nonzero; then the congruence $r^4\equiv s^4
\pmod{m}$ is soluble in $r$, with $F(m)$ solutions $(\mbox{mod }m)$,
where $F$ is a multiplicative function with
\[
F(p)=\left\{
\begin{array}{l l}
1 & \mbox{if $p=2$}
\\
2 & \mbox{if $p\equiv 3\hspace{0.05in}(\mbox{mod }4)$}
\\
4 & \mbox{if $p\equiv 1\hspace{0.05in}(\mbox{mod }4)$}
\end{array}
\right..
\]
Given a solution $r\equiv c\pmod{m}$, we define
\[
N_{c,m,s}=\#\left\{ r: \mbox{$R_i<r\le2R_i$, $\gcd(r,s)=1$ and
$r\equiv c\hspace{-0.1in}\pmod{m}$}\right\}
\]
and
\[
N_{c,s}(n)=\#\left\{ r: \mbox{$R_i<r\le2R_i$, $n|r$ and $r\equiv
c\hspace{-0.1in}\pmod{m}$}\right\};
\]
then
\[
N_{c,m,s}=\sum_{n|s}\mu(n)N_{c,s}(n).
\]
Let $n|s$. Then $\gcd(n,m)=1$, since $\gcd(m,s)=1$. Thus
\[
N_{c,s}(n)=\frac{R_i}{nm}+O(1).
\]
The contribution to $N_{c,m,s}$ of the error term above is of order
at most $d(s)$; that is,
\[
N_{c,m,s}=\frac{R_i}{m}\sum_{n|s}\frac{\mu(n)}{n}+O\left(d(s)\right)=\frac{R_if(s)}{m}+O\left(d(s)\right)
\gg\frac{R_if(s)}{m}.
\]
(The above bound follows from the fact that $m\le R_i^{1/2}$.) Thus
\[
N_{m,s}\gg\frac{F(m)R_if(s)}{m},
\]
and
\[
\displaystyle\sum_{\begin{subarray}{c}R_i<r\le2R_i\\\gcd(r,s)=1\end{subarray}}d\left(\Delta_{r,s}\right)f(\Delta_{r,s})^2
\gg R_if(s)\hspace{-0.1in}\displaystyle\sum_{\begin{subarray}{c}m\le
R_i^{1/2}\\\gcd(m,s)=1\end{subarray}}\frac{F(m)f'(m)}{m}.
\]
In view of the bound (\ref{3LB/eq N_UB 1}) and the fact that $R_i$
and $S_i$ are of the same order, we conclude that
\begin{equation}\label{3LB/eq N_UB 2}
N_{U}(B) \gg
B\displaystyle\sum_{i}\frac{1}{S_i}\displaystyle\sum_{\begin{subarray}{c}m\le
R_i^{1/2}\\\gcd(m,6)=1\end{subarray}}\frac{F(m)f'(m)}{m}\displaystyle\sum_{\begin{subarray}{c}S_i<s\le2S_i\\6|s\\\gcd(m,s)=1\end{subarray}}f(s)^2,
\end{equation}
where the $i$ are summed over those values such that the sets $P_i$
are nonempty. (We impose the condition $\gcd(m,6)=1$ for
convenience.)

\medskip\medskip

We proceed to estimate the inner sum on the right-hand side of
(\ref{3LB/eq N_UB 2}). Let $s=6t$ and $S_i/6=T_i$. Then
\[
\sum_{\begin{subarray}{c}S_i<s\le2S_i\\6|s\\\gcd(m,s)=1\end{subarray}}f(s)^2\gg
\sum_{\begin{subarray}{c}T_i<t\le2T_i\\\gcd(m,t)=1\end{subarray}}f(t)^2.
\]
By the Cauchy-Schwarz inequality, we have
\[
\displaystyle\sum_{\begin{subarray}{c}T_i<t\le2T_i\\\gcd(m,t)=1\end{subarray}}f(t)^2
\ge
\left(\displaystyle\sum_{\begin{subarray}{c}T_i<t\le2T_i\\\gcd(m,t)=1\end{subarray}}
1\right)^{-1}\left(\displaystyle\sum_{\begin{subarray}{c}T_i<t\le2T_i\\\gcd(m,t)=1\end{subarray}}
f(t)\right)^2.
\]
We estimate the two sums on the right-hand side, using the following
two standard relations: first, given any positive integer constant
$c$, we have
\begin{equation}\label{3LB/eq n<N,(n,c)=1}
\#\left\{ n : \mbox{$n\le N$ and
$\gcd(n,c)=1$}\right\}=\frac{N\phi(c)}{c}+O\left(c^\varepsilon\right)
\end{equation}
for any $\varepsilon>0$; and second,
\begin{equation}\label{3LB/lem phi average order}
\displaystyle\sum_{n\le N}\phi(n)=\frac{3N^2}{\pi^2}+O\left(N\log
N\right).
\end{equation}
For the first sum on the right-hand side of the Cauchy-Schwarz
inequality, we have, by (\ref{3LB/eq n<N,(n,c)=1}),
\[
\displaystyle\sum_{\begin{subarray}{c}T_i<t\le2T_i\\\gcd(m,t)=1\end{subarray}}
1 \ll\frac{T_i\phi(m)}{m} = T_if(m).
\]
For the second sum, we have
\[
\displaystyle\sum_{\begin{subarray}{c}T_i<t\le2T_i\\\gcd(m,t)=1\end{subarray}}
f(t)=\displaystyle\sum_{\begin{subarray}{c}T_i<t\le2T_i\\\gcd(m,t)=1\end{subarray}}
\frac{\phi(t)}{t}\gg\frac{1}{T_i}\displaystyle\sum_{\begin{subarray}{c}T_i<t\le2T_i\\\gcd(m,t)=1\end{subarray}}
\phi(t).
\]
By (\ref{3LB/eq n<N,(n,c)=1}), we have
\[
\displaystyle\sum_{\begin{subarray}{c}T_i<t\le2T_i\\\gcd(m,t)=1\end{subarray}}
\phi(t) \ge
\displaystyle\sum_{\begin{subarray}{c}T_i<t\le2T_i\\\gcd(ms,t)=1\end{subarray}}
\phi(t) \gg T_i\displaystyle\sum_{s\le T_i}
\displaystyle\frac{\phi(ms)}{ms} \ge T_if(m)\displaystyle\sum_{s\le
T_i} f(s);
\]
that is,
\[
\displaystyle\sum_{\begin{subarray}{c}T_i<t\le2T_i\\\gcd(m,t)=1\end{subarray}}
f(t)\gg f(m)\displaystyle\sum_{s\le T_i} f(s)\gg T_if(m),
\]
where the second inequality follows from (\ref{3LB/lem phi average
order}). Thus
\[
\displaystyle\sum_{\begin{subarray}{c}T_i<t\le2T_i\\\gcd(m,t)=1\end{subarray}}f(t)^2
\gg T_if(m)\gg S_if(m);
\]
and, in view of (\ref{3LB/eq N_UB 2}),
\begin{equation}\label{3LB/eq N_UB 3}
N_{U}(B) \gg
B\displaystyle\sum_{i}\displaystyle\sum_{\begin{subarray}{c}m\le
R_i^{1/2}\\\gcd(m,6)=1\end{subarray}}\frac{F(m)f(m)f'(m)}{m},
\end{equation}
where the $i$ are summed over those values such that the sets $P_i$
are nonempty.

\medskip\medskip

We now estimate the inner sum on the right-hand side of (\ref{3LB/eq
N_UB 3}). Since $F$, $f$ and $f'$ are all multiplicative, we
consider the corresponding Dirichlet series
\[
D(z)=\displaystyle\sum_{\begin{subarray}{c}m\ge1\\\gcd(m,6)=1\end{subarray}}\frac{F(m)f(m)f'(m)}{m^z},
\]
which admits an Euler product
\[
D(z)=\prod_{p\ge5}\left(1+\displaystyle\frac{F(p)f(p)f'(p)}{p^z}+\sum_{e\ge2}\frac{F(p^e)f(p^e)f'(p^e)}{p^{ez}}\right),
\]
where the product is taken over all primes $p\ge 5$. It is
straightforward to rewrite this as $D(z)=\zeta(z)^3L(z,\chi)F'(z)$,
where $F'$ is a holomorphic function bounded on the half-plane
$\mbox{Re}(z)>3/4$. Hence, by Perron's formula, the inner sum on the
right-hand side of (\ref{3LB/eq N_UB 3}) is equal to
\[
\displaystyle\frac{1}{2\pi
i}\displaystyle\int_{\varepsilon-iT}^{\varepsilon+iT}\zeta(1+w)^3
L(1+w,\chi)F'(1+w)\frac{M^w}{w}dw + O(1).
\]
The integrand has a pole of order $4$ at $w=0$. We apply the residue
theorem to the rectangular contour with corners at $\varepsilon-iT$,
$\varepsilon+iT$, $-1/8+iT$ and $-1/8-iT$, and use the bounds
\[
\zeta(w), L(w,\chi)\ll |w|^{1/8},
\]
which hold provided $\mbox{Re}(w)\ge 7/8$ and $|w-1|\ge 1/8$. These
bounds imply that the integrand along the horizontal segments is of
order at most
\[
\left(T^{1/8}\right)^3T^{1/8}\displaystyle\frac{M^{\mbox{\scriptsize{Re}}(w)}}{T},
\]
where $-1/8\le\mbox{Re}(w)\le\varepsilon$; that is, the contribution
of the integral along the horizontal segments of our contour is of
order at most $1$. Similarly, the integrand along the vertical
segment joining $-1/8+iT$ to $-1/8-iT$ is of order at most
\[
\displaystyle\frac{\left(T^{1/8}\right)^3T^{1/8}}{M^{1/8}};
\]
that is, the contribution of the integral along that segment is of
order at most
\[
\displaystyle\frac{T^{3/2}}{M^{1/8}}=M^{3\varepsilon-1/8}\ll1
\]
provided $\varepsilon<1/24$. Hence we have
\[
\displaystyle\sum_{\begin{subarray}{c}m\le
M\\\gcd(m,6)=1\end{subarray}}\frac{F(m)f(m)f'(m)}{m} \gg \left(\log
M\right)^3 = \left(\log R_i^{1/2}\right)^3.
\]
We insert the above bound into (\ref{3LB/eq N_UB 3}), getting
\[
N_{U}(B) \gg B\displaystyle\sum_{i}\left(\log R_i^{1/2}\right)^3,
\]
where the $i$ are summed over those values such that the sets $P_i$
are nonempty. Now a set $P_i$ is nonempty provided $2^{i+1}\le
B^{1/100}$; that is, provided we have $i\le k\log B$ for some fixed
constant $k>0$. Thus we have:
\[
N_{U}(B) \gg B\displaystyle\sum_{i\le k\log B}\left(\log
R_i^{1/2}\right)^3\gg B\displaystyle\sum_{i\le k\log
B}\left(i-1\right)^3 \gg B\left(\log B\right)^4.
\]

\medskip\medskip

\section{The upper bound}

\medskip

\subsection{Preliminaries}

\medskip\medskip

We define the following projections from $V$ onto $\mathbb{P}^{1}$:
\[
f^{(1)}: [x_{1},\ldots,x_{5}] \mapsto \left\{
           \begin{array}{l l}
              [x_{1},x_{3}]&\mbox{if $(x_{1},x_{3})\neq (0,0)$}\\
              \mbox{$[x_{4},x_{2}]$}&\mbox{otherwise}
           \end{array}
         \right.,
\]
\[
f^{(2)}: [x_{1},\ldots,x_{5}] \mapsto \left\{
           \begin{array}{l l}
              [x_{1},x_{4}]&\mbox{if $(x_{1},x_{4})\neq (0,0)$}\\
              \mbox{$[x_{3},x_{2}]$}&\mbox{otherwise}
           \end{array}
         \right..
\]

\medskip

\noindent We have the following lemma:

\begin{lemma}\label{4UB/lem x height under projections}
$\|f^{(1)}(\textbf{\emph{x}})\| \cdot \|f^{(2)}(\textbf{\emph{x}})\|
\le \|\textbf{\emph{x}}\|$ for all $\textbf{\emph{x}}\in V$.
\end{lemma}

\noindent PROOF. Let $\gcd(x_{1},x_{2},x_{3},x_{4})=n$, and let
$m_{ij}$ denote $\gcd(x_{i},x_{j})n^{-1}$ for $1\le i,j\le4$. Then
\[
\frac{x_1}{n}=m_{13}m_{14} ,\hspace{0.1in}
\frac{x_2}{n}=m_{23}m_{24} ,\hspace{0.1in}
\frac{x_3}{n}=m_{13}m_{23} \hspace{0.1in}\hbox{and}\hspace{0.1in}
\frac{x_4}{n}=m_{14}m_{24}.
\]
Now either $\|f^{(1)}(\textbf{x})\|=\|[x_{1},x_{3}]\|$ or
$\|f^{(1)}(\textbf{x})\|=\|[x_{4},x_{2}]\|$; in both cases we get
$\|f^{(1)}(\textbf{x})\|=\|[m_{14},m_{23}]\|$. Similarly,
$\|f^{(2)}(\textbf{x})\|=\|[m_{13},m_{24}]\|$.\qed

\medskip\medskip

\noindent We define, for $i=1,2$,
\[
N^{(i)}_{U}(B)=\#\{ \textbf{x}\in U: \mbox{$\|\textbf{x}\|\le B$ and
$\|f^{(i)}(\textbf{x})\|\le B^{1/2}$}\}.
\]
\noindent By Lemma \ref{4UB/lem x height under projections}, we have
\[
N_{U}(B)\le N^{(1)}_{U}(B)+N^{(2)}_{U}(B).
\]
We will bound the $N^{(i)}_{U}(B)$. Indeed it suffices to bound
$N^{(1)}_{U}(B)$; the bound for $N^{(2)}_{U}(B)$ follows by
symmetry.

\medskip\medskip

Suppose \textbf{x} contributes to $N^{(1)}_{U}(B)$; say
$f^{(1)}(\textbf{x})=[r,s]$ with $r$ and $s$ coprime. Then
$\textbf{x}$ is of the form $[rX_{1},sX_{2},sX_{1},rX_{2},x_{5}]$,
where $X_{1}=\gcd(x_{1},x_{3})$ and $X_{2}=\gcd(x_{2},x_{4})$; and,
setting $x_{5}=X_{3}$, the second quadric of $V$ is a ternary
quadric $0=Q_{r,s}^{(1)}(\textbf{X})$, where
\[
Q^{(1)}_{r,s}(\textbf{X})=(r^{2}+s^{2})X_{1}^{2}-(r^{2}-s^{2})X_{2}^{2}-2X_{3}^{2}.
\]
The condition $\|\textbf{x}\|\le B$ implies
\begin{equation}\label{4UB/eq X height}
|X_{1}|,|X_{2}|\leq \frac{B}{\mbox{max}(r,s)}
\hspace{0.1in}\mbox{and}\hspace{0.1in} |X_{3}|\leq B.
\end{equation}
Thus, defining
\[
N_{r,s} = \# \{ \mbox{\textbf{X} : $\gcd\left(X_1,X_2,X_3\right)=1$,
$0=Q^{(1)}_{r,s}(\textbf{X})$ and (\ref{4UB/eq X height}) holds} \},
\]
we have
\[
N^{(1)}_{U}(B) \le \displaystyle\sum_{\begin{subarray}{c}
\gcd(r,s) =1 \\
1\le r,s \leq B^{1/2} \end{subarray}} N_{r,s}.
\]

\medskip\medskip

We split the set of suitable pairs $(r,s)$ into dyadic ranges,
letting $P_{i,j}$ denote the set of coprime pairs $(r,s)$ in the
range
\[
\hspace{0.1in} 2^{i-1}=R_{i}<r\le 2R_{i} = 2^{i}
\hspace{0.1in}\hbox{and}\hspace{0.1in} 2^{j-1}=S_{j}<s\le 2S_{j} =
2^{j}.
\]
The bounds $1\le r,s\le B^{1/2}$ imply that the indices $i$ and $j$
have an upper bound $i,j\le k\log B$ for some fixed constant $k>0$.
Thus we have
\begin{equation}\label{4UB/eq comp N_1, N_rs}
N^{(1)}_{U}(B) \ll \displaystyle\sum_{i\le k\log B}
\displaystyle\sum_{j\le i} \sum_{(r,s)\in P_{i,j}}
N_{r,s}+\displaystyle\sum_{j\le k\log B} \displaystyle\sum_{i\le j}
\sum_{(r,s)\in P_{i,j}} N_{r,s}.
\end{equation}
We bound the first of the terms on the right-hand side; the second
term is dealt with similarly.

\medskip

\subsection{Tools}

\medskip\medskip

Our first tool, used to estimate $N_{r,s}$, may be found in
[\cite{Browning/Heath-Brown 2005}]:

\begin{theorem} \label{4UB/thm Heath-Brown}
Let $f \in \mathbb{Z}[\emph{\textbf{X}}]$ be a ternary quadratic
form. Let $M$ denote its matrix representation $M$, $\Delta = |\det
M|\neq0$, and $\Delta_{0}$ denote the highest common factor of the
$2 \times 2$ minors of $M$. Let
\[
N=\#\left\{ \mbox{\emph{$\textbf{X}$ :
$\gcd\left(X_1,X_2,X_3\right)=1$, $0=f(\textbf{x})$ and $|x_{i}|\leq
B_{i}$ for $i=1,2,3$}} \right\}.
\]
Then
\[
N \ll \left( 1 + \left( \frac{B_{1}B_{2}B_{3}\Delta_{0}^{2}}{\Delta}
\right)^{1/3} \right) d(\Delta).
\]
\end{theorem}

We require some notation for our next result. Given
$f\in\mathbb{Z}[x]$ with no fixed prime divisors, the multiplicative
function $\rho_f(m)$ denotes the number of solutions $n\pmod m$ of
$f(n)\equiv0\pmod m$. We collect here some useful results on this
function. The first three are classical, and may be found in
[\cite{Nagell 1981}], for example. The last is attributed in
[\cite{de la Breteche/Browning 2006 III}] to unpublished work by
Stephan Daniel.

\begin{lemma} \label{4UB/lem rho properties}
Let $f\in\mathbb{Z}[x]$ be of degree $g$, have no fixed prime
divisors, and be such that $\mbox{\emph{Disc}}(f)\neq0$. Then:
\end{lemma}

\begin{enumerate}
\item[(a)]
\textit{$\rho_{f}(p)\le g$;}
\item[(b)]
\textit{$\rho_{f}(p^{e})\le gp^{e-1}$ for all $e\in\mathbb{N}$;}
\item[(c)]
\textit{$\rho_{f}(p^{e})=\rho_{f}(p)$ for all $e\in\mathbb{N}$,
provided $p\nmid \mbox{\emph{Disc}}(f)$; and}
\item[(d)]
\textit{$\rho_{f}(p^{e})\le2g^3p^{e\left(1-1/g\right)}$ for all
$e\in\mathbb{N}$.}
\end{enumerate}

\medskip\medskip

\noindent We are now ready to prove the following:

\begin{theorem} \label{4UB/thm Nair}
Let $f\in \mathbb{Z}[x]$ be of degree $4$, have no fixed prime
divisors, and be such that $\mbox{\emph{Disc}}(f)\neq0$. Let
$\alpha, \beta \in (0,1)$ and $N_1,N_2\ge 2$ be such that
$N_2^\alpha\le N_2-N_1\le N_2$ and $\|f\|^\beta \le N_2$. Then the
sum
\[
\displaystyle\sum_{N_1<n\le N_2} d(|f(n)|)
\]
is of order at most
\[
\left(N_2-N_1\right) \displaystyle\prod_{p \le N_2} \left( 1 -
\frac{\rho_{f}(p)}{p} \right) \exp \left( \sum_{p\le N_2}
\frac{d(p)\rho_{f}(p)}{p}
+c\sum_{p|\mbox{\scriptsize{\emph{Disc}}}(f)} \frac{1}{p}\right)
\]
for a constant $c>0$, where the implied constant depends only on
$\alpha$ and $\beta$.
\end{theorem}

\noindent PROOF. This is a special case of the main theorem in
[\cite{Nair 1992}]. Nair's bound depends implicitly on the
discriminant $\mbox{Disc}(f)$. This dependence arises in two places
in [\cite{Nair 1992}]. In both instances we may make explicit or
remove this dependence.

\medskip\medskip

The first instance is in [\cite{Nair 1992}, Lemma 2], in the implied
constant of the bound
\begin{equation}\label{4UB/eq Nair Disc(f) dependence}
\displaystyle\sum_{n\le
N}\left(\frac{n}{\phi(n)}\right)^{4}\frac{d(n)\rho_{f}(n)}{n} \ll
\exp \left( \sum _{p\le
N}\left(\frac{p}{\phi(p)}\right)^{4}\frac{d(p)\rho_{f}(p)}{p}
\right).
\end{equation}
We make this dependence explicit. We begin with the fact that
\[
\displaystyle\sum_{n\le
N}\left(\frac{n}{\phi(n)}\right)^{4}\frac{d(n)\rho_{f}(n)}{n} \le
\exp \left(\displaystyle\sum_{p\le N}
\left(\frac{p}{\phi(p)}\right)^{4}\displaystyle\sum_{e\ge 1}
\frac{d(p^{e})\rho_{f}(p^{e})}{p^{e}} \right).
\]
We shall make use of the bound
\begin{equation}\label{4UB/eq series identity}
\displaystyle\sum_{e\ge E}\frac{e+1}{n^e} \ll
\left(\displaystyle\frac{1}{n^E}\right)\left(\displaystyle\frac{n}{n-1}\right)^2,
\end{equation}
which holds for all $n\in\mathbb{N}$. (Here the relation $\ll$
depends only on $E$.) Now given $p$ such that
$p\nmid\mbox{Disc}(f)$, by Lemma \ref{4UB/lem rho properties}(a),
Lemma \ref{4UB/lem rho properties}(c) and (\ref{4UB/eq series
identity}), we have
\[
\displaystyle\sum_{e\ge 2} \frac{d(p^{e})\rho_{f}(p^{e})}{p^{e}} \le
4\displaystyle\sum_{e\ge 2} \frac{(e+1)}{p^{e}} \ll \frac{1}{p^2}.
\]
Likewise, given $p$ such that $p|\mbox{Disc}(f)$, by Lemma
\ref{4UB/lem rho properties}(d) and (\ref{4UB/eq series identity}),
we have
\[
\displaystyle\sum_{e\ge8} \frac{d(p^{e})\rho_{f}(p^{e})}{p^{e}}
\le128\displaystyle\sum_{e\ge8} \frac{e+1}{p^{e/4}} \ll
\frac{1}{p^2}.
\]
Finally, given $p$ such that $p|\mbox{Disc}(f)$, by Lemma
\ref{4UB/lem rho properties}(b), we have
\[
\displaystyle\sum_{2\le e<8} \frac{d(p^{e})\rho_{f}(p^{e})}{p^{e}}
\le 4\displaystyle\sum_{e<8} \frac{e+1}{p} \ll \frac{1}{p}.
\]
These bounds combine to give
\[
\displaystyle\sum_{n\le
N}\left(\frac{n}{\phi(n)}\right)^4\frac{d(n)\rho_{f}(n)}{n} \ll \exp
\left( \sum _{p\le
N}\left(\frac{p}{\phi(p)}\right)^4\frac{d(p)\rho_{f}(p)}{p} +
c\sum_{p|\mbox{\scriptsize{Disc}}(f)} \frac{1}{p} \right)
\]
for a constant $c>0$, where the relation $\ll$ does not depend on
$\mbox{Disc}(f)$. The difference between this bound and (\ref{4UB/eq
Nair Disc(f) dependence}) accounts for the difference between
Theorem \ref{4UB/thm Nair} and the main result in [\cite{Nair
1992}].

\medskip\medskip

The second place in [\cite{Nair 1992}] in which a dependence on
$\mbox{Disc}(f)$ arises is in the author's reduction of the bound
[\cite{Nair 1992}, (6.3)], where, given a positive integer $n$ such
that $N^{1/2}<n\le N$, the bound $\rho_{f}(n)\ll N^{1/8}$ is
invoked; $\mbox{Disc}(f)$ figures in the implied constant. We remove
the dependence on $\mbox{Disc}(f)$ by invoking the bound
$\rho_{f}(n)\ll n^{4/5}$ for all $n\in\mathbb{N}$, where the
relation $\ll$ does not depend on $\mbox{Disc}(f)$; this proves to
be sufficient. \qed

\medskip\medskip

\noindent We use Theorem \ref{4UB/thm Nair} to prove our version of
a result due to Browning and de la Bret\`{e}che, which we use to sum
our estimates for $N_{r,s}$ over the pairs $(r,s)\in P_{i,j}$. We
require a generalization of the function $\rho_f$ to binary forms.
Let $f\in\mathbb{Z}[x_1,x_2]$ have no fixed prime divisors. Then
$\rho_{f(1,x)}(m)$ denotes the number of solutions $n\pmod m$ of
$f(1,n)\equiv0\pmod m$, and we define for any prime $p$ the function
\[
\rho^{*}_{f}(p)=\left\{ \begin{array}{l l} \rho_{f(1,x)}(p)+1 &
\mbox{if $p|f(0,1) $}
\\ \rho_{f(1,x)}(p) & \mbox{otherwise} \end{array}
\right..
\]

\begin{theorem} \label{4UB/thm Breteche/Browning}
Let $f\in \mathbb{Z}[x_{1},x_{2}]$ be of degree $4$, have no fixed
prime divisors, and be such that $\mbox{\emph{Disc}}(f)\neq 0$ and
$f(1,0)f(0,1)\neq0$. Let $\alpha,\beta\in(0,1)$ and $N,N_1,N_2\ge2$
be such that $N_2^{\alpha}\le N_2-N_1\le N_2$ and $\min(N,N_2)\ge
a\max(N,N_2)^{4\beta}\|f\|^{\beta}$ for a constant $a>0$ dependent
only on $\beta$. Then
\[
\displaystyle\sum_{1\le n_1\le N} \sum_{N_1< n_2\le N_2}
d\left(|f(n_1,n_2)|\right) \ll N(N_2-N_1)T,
\]
where
\[
T=\displaystyle\prod_{p|\mbox{\scriptsize{\emph{Disc}}}(f)}\left(1+\frac{1}{p}\right)^{b}
\exp\left(c\sum_{p|\mbox{\scriptsize{\emph{Disc}}}(f)}
\frac{1}{p}\right)\exp\left(\sum_{p\le\max(N,N_2)}\frac{\rho^{*}_{f}(p)}{p}\right)
\]
for constants $b,c>0$, and the relation $\ll$ depends only on
$\alpha$ and $\beta$.
\end{theorem}

\noindent PROOF. This theorem is an adaptation of [\cite{de la
Breteche/Browning 2006 III}, Theorem 1]. There, the authors take
$n_2\le N_2$; we take a shorter range of summation and appeal to
Theorem \ref{4UB/thm Nair}. As in [\cite{de la Breteche/Browning
2006 III}, \S3], we fix $n_1$ and consider the sum
\[
\sum_{N_1< n_{2}\le N_2} d\left(|f(n_{2})|\right).
\]
By Theorem \ref{4UB/thm Nair}, the sum above has an upper bound of
order at most
\[
(N_2-N_1)\displaystyle\prod_{p \le N_2} \left( 1 -
\frac{\rho_{f}(p)}{p} \right)\exp\left(\displaystyle\sum_{p\le N_2}
\frac{d(p)\rho_{f}(p)}{p} +c\sum_{p|n_1\mbox{\scriptsize{Disc}}(f)}
\frac{1}{p}\right)
\]
for a constant $c>0$. In comparison, in [\cite{de la
Breteche/Browning 2006 III}, \S3] the authors conclude that
\[
\sum_{n_{2}\le N_2} d\left(|f(n_{2})|\right)\ll
N_2\displaystyle\prod_{p \le N_2} \left( 1 - \frac{\rho_{f}(p)}{p}
\right)\displaystyle\sum_{n_2\le N_2}\frac{d(n_2)\rho_f(n_2)}{n_2}.
\]
This difference accounts for the discrepancy between Theorem
\ref{4UB/thm Breteche/Browning} and [\cite{de la Breteche/Browning
2006 III}, Theorem 1]. Proceeding according to the argument of
[\cite{de la Breteche/Browning 2006 III}, \S3], we have
\[
\displaystyle\sum_{1\le n_1\le N} \sum_{N_1< n_2\le N_2}
d\left(|f(n_1,n_2)|\right) \ll N(N_2-N_1)t_1t_2,
\]
where
\[
\begin{array}{l}
t_1 =\displaystyle\prod_{4<p\le
N_2}\left(1-\frac{\rho_{f(1,x)}(p)}{p}\right)\exp\left(\sum_{\begin{subarray}{c}p\le
N_2\\p\nmid n_1\end{subarray}}
\frac{d(p)\rho_{f}(p)}{p}\right)\hspace{0.1in}\mbox{and}
\smallskip
\\
t_2 =
\displaystyle\prod_{p|\mbox{\scriptsize{Disc}}(f)}\left(1+\frac{1}{p}\right)^{b}
\exp\left(c\sum_{p|\mbox{\scriptsize{Disc}}(f)} \frac{1}{p}\right),
\end{array}
\]
for constants $b,c>0$. It is straightforward to show that $t_1$ is
of order at most
\[
\displaystyle\prod_{4<p\le
N_2}\hspace{-0.05in}\left(1-\frac{\rho_{f(1,x)}(p)}{p}\right)\exp\left(\sum_{4<p\le
N_2}
\hspace{-0.05in}\frac{2\rho_{f(1,x)}(p)}{p}\right)\ll\exp\left(\displaystyle\sum_{p\le
N_2} \frac{\rho_{f(1,x)}(p)}{p}\right)
\]
which, combined with $t_2$, yields the theorem. \qed

\medskip\medskip

Our third main tool is a classical result due to Dedekind and
Landau:

\begin{theorem} \label{4UB/thm Dedekind/Landau}
Let $f \in \mathbb{Z}[x]$ be irreducible and of degree $g\ge 1$.
Then
\[
\sum_{p\le B}\rho_{f}(p)=\mbox{\emph{Li}}(B)+O\left(
\displaystyle\frac{B}{\exp(c (\log B)^{1/2})} \right)
\]
for a constant $c>0$ dependent only on the splitting field of $f$
over $\mathbb{Q}$.
\end{theorem}

\noindent PROOF. Let $L$ be the splitting field of $f$ over
$\mathbb{Q}$. For all but finitely many $p$, $f\pmod p$ has
factorization $F_{1}\cdots F_{n}\pmod p$, where the $F_{i}\in
\mathbb{Z}_{p}[x]$ are irreducible and of degrees $g_i$,
respectively, if and only if the principal ideal $(p)$ has
factorization $P_{1}\cdots P_{n}$, where the $P_{i}$ are prime
ideals over $L$ with norms $p^{g_i}$, respectively. Now
\[
\rho_f(p)= \#\{i: F_{i}\mbox{ is linear} \} = \#\{i:
\mbox{norm}(P_{i})=p \},
\]
and by Landau's Prime Ideal Theorem,
\[
\#\{ \mbox{prime ideals } P : \mbox{norm}(P)=p\le B  \} =
\mbox{Li}(B)+O\left( \displaystyle\frac{B}{\exp(c (\log B)^{1/2})}
\right)
\]
for a constant $c>0$ dependent only on $L$. \qed

\medskip\medskip

\noindent Now the Prime Ideal Theorem is simply the generalization
to number fields of the Prime Number Theorem; given
$n\in\mathbb{N}$, we have
\[
\pi(n) = \mbox{Li}(n) + O\left( \frac{n}{\exp(c'(\log n)^{1/2})}
\right)
\]
for a constant $c'>0$. This symmetry between $\pi(t)$ and the
average order of $\rho_{f}(p)$ will be useful. We also record the
following bound, due to Rosser and Schoenfeld
[\cite{Rosser/Schoenfeld 1962}]:

\begin{lemma} \label{4UB/lem Rosser/Schoenfeld} Let $n \ge 67$.
Then $\displaystyle\frac{n}{\log n - 1/2} < \pi(n) < \frac{n}{\log n
- 3/2}$.
\end{lemma}

\medskip

\subsection{The proof of the upper bound}

\medskip\medskip

As in \S3, we let $\Delta_{r,s}$ denote $|r^4-s^4|$. We shall also
write $P$, $R$ and $S$ for $P_{i,j}$, $R_i$ and $S_j$, respectively.

\medskip\medskip

We begin by applying Theorem \ref{4UB/thm Heath-Brown}, getting
\begin{equation}\label{4UB/eq Heath-Brown applied 2}
\displaystyle\sum_{(r,s)\in P} N_{r,s} \ll B \left(
\displaystyle\frac{1}{R^{2/3}} \displaystyle\sum_{(r,s)\in P}
\displaystyle\frac{d(\Delta_{r,s})}{\Delta_{r,s}^{1/3}}\right).
\end{equation}
We evaluate the sum on the right-hand side according to the size of
$\Delta_{r,s}$. Let the linear factors of $\Delta_{r,s}$ be denoted
$|s-\alpha_ir|$ for $i=1,2,3,4$. We consider three cases:
\[
\begin{array}{l l}
\mbox{Case I:} & \mbox{$R$ and $S$ are not of the same order;}
\\
\mbox{Case II:} & \mbox{$R$ and $S$ are of the same order, and
$|s-\alpha_ir|>R/4$}
\\
& \mbox{for $i=1,2,3,4$; and}
\\
\mbox{Case III:} & \mbox{$R$ and $S$ are of the same order, and
$|s-\alpha_ir|\le R/4$}
\\
& \mbox{for some $i\in\{1,2,3,4\}$. (We may assume moreover that}
\\
& \mbox{$\alpha_i=1$ or $-1$, for otherwise we have
$|s-\alpha_ir|>R/4$.)}
\end{array}
\]
Note that, since we are in search of an upper bound, we may apply
selectively the coprimality condition on $P$.

\medskip\medskip

In Case I, $\Delta_{r,s}$ is dominated by the $r^4$ term, and we
have
\[
\displaystyle\sum_{(r,s)\in P} N_{r,s} \ll B \left(
\displaystyle\frac{1}{R^2} \displaystyle\sum_{(r,s)\in P}
d(\Delta_{r,s})\right).
\]
Now
\[
\displaystyle\sum_{(r,s)\in P} d(\Delta_{r,s})\le
\displaystyle\sum_{s\le
2\max\left(S,R^{1/2}\right)}\displaystyle\sum_{r\le 2R}
d(\Delta_{r,s}).
\]
We apply Theorem \ref{4UB/thm Breteche/Browning} to the right-hand
side, getting
\[
\displaystyle\sum_{(r,s)\in P} d(\Delta_{r,s})\ll
\max\left(S,R^{1/2}\right)RT,
\]
where
\[
T=\displaystyle\prod_{p|\mbox{\scriptsize{Disc}}(\Delta_{r,s})}\left(1+\frac{1}{p}\right)^{b}
\exp\left(c\sum_{p|\mbox{\scriptsize{Disc}}(\Delta_{r,s})}
\frac{1}{p}\right)\exp\left(\sum_{p\le
2R}\frac{\rho^{*}_{\Delta_{r,s}}(p)}{p}\right)
\]
for some constants $b,c>0$. The fact that
$\mbox{Disc}\left(\Delta_{r,s}\right)=128i$ implies that the first
two terms are $\ll1$, whence
\[
T \ll \exp\left(\sum_{p\le
2R}\frac{\rho^{*}_{\Delta_{r,s}}(p)}{p}\right)\ll
\exp\left(\displaystyle\sum_{p\le2R}
\frac{\rho_{\Delta_{1,x}}(p)}{p}\right).
\]
We appeal to Theorem \ref{4UB/thm Dedekind/Landau} and Lemma
\ref{4UB/lem Rosser/Schoenfeld}. The sum on the far right-hand side
is equal to
\[
\frac{1}{2R}\sum_{p\le2R}\rho_{\Delta_{1,x}}(p) +
\int_{1}^{2R}\sum_{p\le t}\rho_{\Delta_{1,x}}(p)\frac{dt}{t^2} +
O\left( \int_{1}^{2R}\sum_{p\le
t}\rho_{\Delta_{1,x}}(p)\frac{dt}{t^3}\right).
\]
The first term is small. Indeed, let $f_1\left(x\right)=1+x^2$,
$f_2\left(x\right)=1+x$ and $f_3\left(x\right)=1-x$. Then, by
Theorem \ref{4UB/thm Dedekind/Landau}, the first term is equal to
\[
\frac{1}{2R}\left(\sum_{p\le2R}\rho_{f_{1}(x)}(p)+\sum_{p\le2R}\rho_{f_{2}(x)}(p)+\sum_{p\le2R}\rho_{f_{3}(x)}(p)\right)
=O\left(1\right).
\]
The error term is also small: by Lemma \ref{4UB/lem rho
properties}(a) we have $\rho_{\Delta_{1,x}}(p)\le 4$ for all primes
$p$; that is,
\[
O\left( \int_{1}^{2R}\sum_{p\le
t}\rho_{\Delta_{1,x}}(p)\frac{dt}{t^3}\right)=O\left(\int_{1}^{2R}\frac{dt}{t^2}\right)=O\left(1\right).
\]
Thus we have
\[
T \ll \exp\left(\displaystyle\int_{67}^{2R}\sum_{p\le
t}\rho_{\Delta_{1,x}}(p)\frac{dt}{t^2}+O\left(1\right)\right).
\]
By Lemma \ref{4UB/lem Rosser/Schoenfeld}, we have
\begin{equation}\label{4UB/eq Case I 2}
\int_{67}^{2R}\sum_{p\le t}\rho_{\Delta_{1,x}}(p)\frac{dt}{t^{2}} <
\int_{67}^{2R}\frac{1}{\pi(t)}\sum_{p\le
t}\rho_{\Delta_{1,x}}(p)\frac{dt}{t\left(\log t -3/2\right)}.
\end{equation}
By Theorem \ref{4UB/thm Dedekind/Landau}, we have
\[
\begin{array}{r c l}
\displaystyle\sum_{p\le t}\rho_{\Delta_{1,x}}(p) & = &
\displaystyle\sum_{p\le t}\rho_{f_{1}(x)}(p) + \sum_{p\le
t}\rho_{f_{2}(x)}(p) + \sum_{p\le t}\rho_{f_{3}(x)}(p)\medskip
\\
& = & 3 \left( \hbox{Li}(t) + O\left( \displaystyle\frac{t}{\exp(c
(\log t)^{1/2})} \right) \right)
\end{array}
\]
for a constant $c>0$. We also have
\[
\pi(t) = \mbox{Li}(t) + O\left( \frac{t}{\exp(c'(\log t)^{1/2})}
\right)
\]
for a constant $c'>0$. Let $C=\min(c,c')$. Then
\[
\displaystyle\frac{1}{\pi(t)} \sum_{p\le t}\rho_{\Delta_{x,1}}(p) =
3+ O\left(\displaystyle\frac{1}{\log t-3/2}\right).
\]
Substituting back into (\ref{4UB/eq Case I 2}), we get
\[
\begin{array}{r c l}
\displaystyle\int_{67}^{2R}\sum_{p\le
t}\rho_{\Delta_{1,x}}(p)\frac{dt}{t^{2}} &<&
\displaystyle\int_{67}^{2R}\frac{3\,dt}{t\left(\log t -3/2\right)} +
O\left(\int_{67}^{2R}\frac{dt}{t\left(\log t -3/2\right)^2}\right)
\smallskip
\\
&=& \log\left(\log(2R)-2/3\right)^3+O(1)
\smallskip
\\
&<& \log\left(\log B\right)^3+O(1).
\end{array}
\]
Thus $T\ll\left(\log B\right)^3$, and
\[
\displaystyle\sum_{(r,s)\in
P}d(\Delta_{r,s})\ll\max\left(S,R^{1/2}\right)R\left(\log
B\right)^3;
\]
that is,
\begin{equation}\label{4UB/eq Case I final}
\displaystyle\sum_{(r,s)\in P} N_{r,s}
\ll\max\left(\frac{S}{R},\frac{1}{R^{1/2}}\right)B\left(\log
B\right)^3
\end{equation}
for Case I.

\medskip\medskip

Case II is handled identically: as in Case I, we have
$\Delta_{r,s}\gg R^{4}$, and the same bound (\ref{4UB/eq Case I
final}) results.

\medskip\medskip

In Case III, suppose $\alpha_1\in\mathbb{R}$ and $|s-\alpha_1r|\le
R/4$. Then the bounds
\[
r|\alpha_{1}-\alpha_{i}|-|s-\alpha_{1}r|
 \le |s-\alpha_{i}r|\le r|\alpha_{1}-\alpha_{i}| + |s-\alpha_{1}r|
\]
for $i=2,3,4$ imply that $\Delta_{r,s}$ is of order
$|s-\alpha_{1}r|R^{3}$. We split the set of values for
$|s-\alpha_{1}r|$ into dyadic ranges
\[
2^{i-1}=B_{i}< |s-\alpha_{1}r| \le 2B_{i} = 2^{i},
\]
where the index $i$ has an upper bound
\[
I=\left\lceil\displaystyle\frac{\log(
R/4)}{\log2}\right\rceil=\frac{\log R}{\log 2}+O(1).
\]
In view of (\ref{4UB/eq Heath-Brown applied 2}), we have
\[
\displaystyle\sum_{(r,s)\in P} N_{r,s} \ll
B\left(\displaystyle\frac{1}{R^{5/3}} \displaystyle\sum_{i\le
I}\frac{1}{B_{i}^{1/3}}\sum_{\begin{subarray}{c}(r,s)\in P
\\ B_{i}<|s-\alpha_{1}r|\le2B_{i} \end{subarray}}d(\Delta_{r,s})\right).
\]
Now
\begin{equation}\label{4UB/eq Case III 1}
\displaystyle\sum_{\begin{subarray}{c}(r,s)\in P
\\ B_{i}<|s-\alpha_{1}r|\le2B_{i} \end{subarray}} d(\Delta_{r,s})
\le\displaystyle\sum_{1\le s\le 2S} \displaystyle\sum_{K_{i}\le r\le
L_{i}} d(\Delta_{r,s})
\end{equation}
where
\[
K_{i}=\max\left(1,s-2\max\left(B_{i},S^{1/3}\right)\right)
\]
and
\[
L_{i}=\min\left(2R,s+2\max\left(B_{i},S^{1/3}\right)\right).
\]
We apply Theorem \ref{4UB/thm Breteche/Browning} to the right-hand
side of (\ref{4UB/eq Case III 1}), getting
\[
\displaystyle\sum_{\begin{subarray}{c}(r,s)\in P
\\ B_{i}<|s-\alpha_{1}r|\le2B_{i} \end{subarray}} d(\Delta_{r,s})\ll\max\left(B_{i},S^{1/3}\right) S \left(\log
B\right)^3,
\]
hence
\[
\displaystyle\sum_{(r,s)\in P} N_{r,s} \ll
\left(\displaystyle\frac{S}{R^{5/3}} \displaystyle\sum_{i\le
I}\frac{\max\left(B_{i},S^{1/3}\right)}{B_{i}^{1/3}}
\right)B\left(\log B\right)^3.
\]

\medskip\medskip

If $\max\left(B_{i},S^{1/3}\right)=B_i$, we have
\[
\displaystyle\sum_{i\le
I}\frac{\max\left(B_{i},S^{1/3}\right)}{B_{i}^{1/3}} \ll
2^{2I/3}=\exp\left(\frac{2I}{3}\log 2\right) \ll R^{2/3};
\]
and if $\max\left(B_{i},S^{1/3}\right)=S^{1/3}$, we have
\[
\displaystyle\sum_{i\ll\log
R}\frac{\max\left(B_{i},S^{1/3}\right)}{B_{i}^{1/3}} \ll S^{1/3}.
\]
Thus we have, for Case III, the bound
\begin{equation}\label{4UB/eq Case III}
\displaystyle\sum_{(r,s)\in P} N_{r,s} \ll
\max\left(\displaystyle\frac{S}{R},\frac{S^{4/3}}{R^{5/3}}\right)
B\left(\log B\right)^3=\left(\frac{S}{R}\right)B\left(\log
B\right)^3.
\end{equation}

\medskip\medskip

Comparing the bounds (\ref{4UB/eq Case I final}) and (\ref{4UB/eq
Case III}), we conclude that
\begin{equation}\label{4UB/eq N_rs}
\displaystyle\sum_{(r,s)\in P}
N_{r,s}\ll\max\left(\displaystyle\frac{S}{R},\frac{1}{R^{1/2}}\right)
B\left(\log B\right)^3.
\end{equation}

\medskip
\section{The cardinality $N_{U}(B)$}\label{4UB/section The cardinality N_WB}
\medskip\medskip

By the bounds (\ref{4UB/eq comp N_1, N_rs}) and (\ref{4UB/eq N_rs}),
we have
\[
\displaystyle\sum_{i\le k\log B} \displaystyle\sum_{j\le i}
\sum_{(r,s)\in P_{i,j}} N_{r,s}\ll B\left(\log
B\right)^3\displaystyle\sum_{i\le k\log B} \displaystyle\sum_{j\le
i}
\max\left(\displaystyle\frac{S_j}{R_i},\frac{1}{R_i^{1/2}}\right).
\]
If $S_j\ge R_i^{1/2}$, the sum on the right-hand side is equal to
\[
\displaystyle\sum_{i\le k\log B}\displaystyle\sum_{j\le
i}\frac{1}{2^{i-j}}\le \displaystyle\sum_{i\le k\log
B}\displaystyle\sum_{j\ge0}\frac{1}{2^j}\ll \log B;
\]
otherwise, it is equal to
\[
\displaystyle\sum_{i\le k\log B}\displaystyle\sum_{j\le
i}\frac{1}{2^{(i-1)/2}}\ll1.
\]
Thus we have
\[
\displaystyle\sum_{i\le k\log B} \displaystyle\sum_{j\le i}
\sum_{(r,s)\in P_{i,j}} N_{r,s}\ll B\left(\log B\right)^4
\]
as required.

\medskip\medskip

\section*{Acknowledgements}

\medskip

The main result in this paper was first established in the author's
doctoral thesis, completed at the University of Oxford under the
supervision of Roger Heath-Brown, and with the financial support of
the Clarendon Fund, the National Sciences and Engineering Research
Council of Canada, and Balliol College.

\medskip\medskip

\noindent \textit{Department of Pure Mathematics}\\
\noindent \textit{Faculty of Mathematics}\\
\noindent \textit{University of Waterloo}\\
\noindent \textit{200 University Avenue West}\\
\noindent \textit{Waterloo, Ontario, Canada N2L 3G1}\\
\noindent \textit{email: fsleung@math.uwaterloo.ca}\\


\begin{thebibliography}{99}\label{BIB}



\bibitem{de la Breteche/Browning 2006 III} de la Bret\`{e}che, R. and
Browning, T.D. ``Sums of arithmetic functions over values of binary
forms.'' \emph{Acta Arith.} \textbf{125} (2006), 291-304.

\bibitem{Browning 2006} Browning, T.D. ``An overview of
Manin's conjecture for del Pezzo surfaces.'' \emph{Proceedings of
the Gauss-Dirichlet Conference}, G\"{o}ttingen, 2006.

\bibitem{Browning/Heath-Brown 2005} Browning, T.D. and Heath-Brown, D.R. ``Counting rational points on
hypersurfaces.'' \emph{J. Reine Angew. Math.} \textbf{584} (2005),
83-115.

\bibitem{Franke/Manin/Tschinkel 1989} J. Franke, Y.I. Manin and Y.
Tschinkel. ``Rational points of bounded height on Fano varieties.''
\emph{Invent. Math.} \textbf{95} (1989), 421-435.

\bibitem{Manin 1974} Manin, Y.I. \emph{Cubic Forms: Algebra, Geometry,
Arithmetic}, North-Holland, Amsterdam, 1974.

\bibitem{Nagell 1981} Nagell, T. \emph{Introduction to Number Theory}, 2nd ed., Chelsea, New York, 1981.

\bibitem{Nair 1992} Nair, M. ``Multiplicative functions of polynomial values in short
intervals.'' \emph{Acta Arith.} \textbf{62} (1992), 257-269.

\bibitem{Rosser/Schoenfeld 1962} Rosser, J.B. and Schoenfeld, L. ``Approximate formulas for some functions of prime
numbers.'' \emph{Illinois J. Math.} \textbf{6} (1962), 64-94.



\end{thebibliography}
\end{document}